# ASYMPTOTIC CAUCHY GAINS: DEFINITIONS AND SMALL-GAIN PRINCIPLE

E.D. Sontag[*]
Department of Mathematics
Rutgers University, New Brunswick, NJ 08903
http://www.math.rutgers.edu/~sontag

**Abstract**

A notion of "asymptotic Cauchy gain" for input/output systems, and an associated small-gain principle, are introduced. A Lyapunov-like characterization allows the computation of these gains for state-space systems, and the formulation of sufficient conditions insuring the lack of oscillations and chaotic behaviors in a wide variety of cascades and feedback loops.

## 1 Introduction

In this note, we introduce a notion of "asymptotic Cauchy gain" for input/output systems, and establish a simple small-gain principle as well as a Lyapunov-like characterization which allows the computation of these gains for state-space systems.

We were motivated by the problem of guaranteeing the non-existence of oscillations in certain biological inhibitory feedback loops. Standard small-gain principles are hard to apply in that context, because the location of closed-loop equilibria may depend on the gain of the feedback law.

References and comparisons to other approaches to small-gain theorems, as well as to other stability results for feedback loops, e.g. based upon the circle criterion and passivity, the Nyquist criterion, or the secant condition, will be discussed in the final version of this preliminary report.

### 1.1 Cauchy Gains

For any metric space $M$, we write the distance $d_M(a,b)$ between any two elements $a, b \in M$, in the suggestive form "$|a-b|$" even when $M$ has no linear structure (so the "$-$" sign has no meaning, of course), and define the *asymptotic amplitude* of a function $\omega : \mathbb{R}_{\geq 0} \to M$, where $\mathbb{R}_{\geq 0} = [0, +\infty)$, as follows:

$$\|\omega\|_{\mathrm{aa}} := \limsup_{s,t \to \infty} |\omega(t) - \omega(s)| = \lim_{T \to \infty} \left( \sup_{t,s \geq T} |\omega(t) - \omega(s)| \right) \in [0, \infty].$$

---

[*]Supported in part by US Air Force Grant F49620-01-1-0063

Observe that the condition "$\|\omega\|_{\mathrm{aa}} = 0$" amounts to the Cauchy property for $\omega$: for every $\varepsilon > 0$ there is some $T > 0$ such that $|\omega(t) - \omega(s)| < \varepsilon$ for all $t, s \geq T$. Thus, when $M$ is a complete metric space (for instance, if, as in all our examples, $M \subseteq \mathbb{R}^m$ is any closed subset of a Euclidean space):

$$\|\omega\|_{\mathrm{aa}} = 0 \quad \Longleftrightarrow \quad \exists \lim_{t \to \infty} \omega(t).$$

If $\|\omega\|_{\mathrm{aa}} = 0$, we denote $\omega^\infty := \lim_{t \to \infty} \omega(t)$.

Let $\mathcal{U}$ and $\mathcal{Y}$ be two complete metric spaces. We define a *behavior with input-value space $\mathcal{U}$ and output-value space $\mathcal{Y}$* as a relation $\mathcal{R}$ between time-functions with values in $\mathcal{U}$ and $\mathcal{Y}$ respectively:

$$\mathcal{R} \subseteq [\mathbb{R}_{\geq 0} \to \mathcal{U}] \times [\mathbb{R}_{\geq 0} \to \mathcal{Y}]$$

where $[\mathbb{R}_{\geq 0} \to M]$ is the set of functions $\mathbb{R}_{\geq 0} \to M$. We call any element $(\omega, \eta) \in \mathcal{R}$ an *input/output pair*, and say that $\omega$ is an input signal and $\eta$ is an output signal of $\mathcal{R}$.

Typical examples of behaviors are those obtained by starting with a system of differential equations with inputs ("forcing functions" or "controls") $\omega$, and viewing the solutions obtained by solving the system with different initial states, or some components of these solutions, as the outputs $\eta$. This will be discussed in detail later.

We use standard terminology for comparison functions: $\mathcal{K}_\infty$ is the class of continuous, strictly increasing, and unbounded functions $\gamma : \mathbb{R}_{\geq 0} \to \mathbb{R}_{\geq 0}$ with $\gamma(0) = 0$.

**Definition 1.1** A behavior $\mathcal{R}$ has *Cauchy gain* $\gamma \in \mathcal{K}_\infty$ if

$$\|\eta\|_{\mathrm{aa}} \leq \gamma(\|\omega\|_{\mathrm{aa}})$$

for all $(\omega, \eta) \in \mathcal{R}$. □

The existence of a Cauchy gain for $\mathcal{R}$ implies, in particular, the following *converging input converging output* property for $\mathcal{R}$: if $\omega(t) \to \bar{u}$ as $t \to \infty$, for some $\bar{u} \in \mathcal{U}$ (that is, if $\|\omega\|_{\mathrm{aa}} = 0$), and if $(\omega, \eta) \in \mathcal{R}$, then also $\eta(t) \to \bar{y}$ as $t \to \infty$, for some $\bar{y} \in \mathcal{Y}$.

The interconnection that results when the output of a system $\mathcal{R}$ is fed back to its input under the action of the system (feedback law) $\mathcal{S}$ is pictorially represented in Figure 1. The behavioral terminology gives an easy way to define formally the meaning of this interconnection: if $\mathcal{R}$ and $\mathcal{S}$ are behaviors, then the signals that appear when the loop is closed are precisely those pairs $(\omega, \eta)$ such that $(\omega, \eta) \in \mathcal{R}$ and $(\eta, \omega) \in \mathcal{S}$. Put another way, the feedback connection is simply the behavior $\mathcal{R} \bigcap \mathcal{S}^{-1}$, where, for any behavior $\mathcal{S}$ with input-value space $\mathcal{Y}$ and output-value space $\mathcal{U}$, we denote by $\mathcal{S}^{-1}$ the inverse behavior, with input-value space $\mathcal{U}$ and output-value space $\mathcal{Y}$, consisting of all pairs $(\omega, \eta)$ such that $(\eta, \omega) \in \mathcal{S}$.

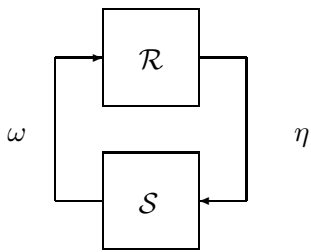

Figure 1: Feedback Interconnection $\mathcal{R} \bigcap \mathcal{S}^{-1}$

With this formalism, the basic "small gain principle" is trivial to establish. It states that the interconnection of two systems having Cauchy gains whose composition is a contraction, has the property that the external signals $\omega$ and $\eta$ must always converge to some value as $t \to \infty$, at least if they are known to have finite asymptotic amplitude:

**Lemma 1.2** *(Small gain lemma for asymptotic amplitude.)* Suppose that $\mathcal{R}$ and $\mathcal{S}$ are two behaviors with Cauchy gains $\gamma_1$ and $\gamma_2$ respectively, and that the following condition holds:

$$\gamma_1(\gamma_2(r)) < r \quad \forall\, r > 0\,. \tag{1}$$

Then, for all $(\omega, \eta) \in \mathcal{R} \bigcap \mathcal{S}^{-1}$ for which $\|\omega\|_{\mathrm{aa}} < \infty$, $\|\omega\|_{\mathrm{aa}} = \|\eta\|_{\mathrm{aa}} = 0$.

*Proof.* Since $(\omega, \eta) \in \mathcal{R}$, $\|\eta\|_{\mathrm{aa}} \leq \gamma_1(\|\omega\|_{\mathrm{aa}})$; and since also $(\eta, \omega) \in \mathcal{S}$, $\|\omega\|_{\mathrm{aa}} \leq \gamma_2(\|\eta\|_{\mathrm{aa}})$. If $\|\eta\|_{\mathrm{aa}} \neq 0$, then $\|\eta\|_{\mathrm{aa}} \leq \gamma_1(\gamma_2(\|\eta\|_{\mathrm{aa}})) < \|\eta\|_{\mathrm{aa}}$, a contradiction. Finally, $\|\omega\|_{\mathrm{aa}} \leq \gamma_2(\|\eta\|_{\mathrm{aa}}) = \gamma_2(0) = 0$ gives that also $\|\omega\|_{\mathrm{aa}} = 0$. ∎

**Remark 1.3** Note that the condition "$\|\omega\|_{\mathrm{aa}} < \infty$" is equivalent to ultimate boundedness, i.e. there are a bounded set $C \subseteq \mathcal{U}$ and some $T \geq 0$ such that $\omega(t) \in C$ for all $t \geq T$. (Writing $|u| := |u - 0|$ for some fixed element $0 \in \mathcal{U}$: if there are some $c, T > 0$ so that $|\omega(t)| \leq c$ for all $t \geq T$ then $\|\omega\|_{\mathrm{aa}} \leq 2c$; conversely, if $\sup_{t,s \geq T} |\omega(t) - \omega(s)| \leq c$ for some $T$ then $|w(t)| \leq c + |\omega(T)|$ for all $t \geq T$.) In applications to feedback loops involving differential equations, all signals are continuous, and for them, ultimate boundedness is equivalent to just boundedness. □

The limiting values of the signals $\omega$ and $\eta$ need not be unique; for instance bistable systems give rise to nonuniqueness. In order to present a condition which guarantees uniqueness, we introduce a new concept.

**Definition 1.4** A behavior $\mathcal{R}$ has *incremental limit gain* $\kappa \in \mathcal{K}_\infty$ if the following property holds:

$$\limsup_{t \to \infty} |\eta_1(t) - \eta_2(t)| \;\leq\; \kappa(|\omega_1^\infty - \omega_2^\infty|)$$

whenever $(\omega_i, \eta_i) \in \mathcal{R}$ are any two pairs with the properties $\|\omega_1\|_{\mathrm{aa}} = \|\omega_2\|_{\mathrm{aa}} = 0$. □

In words, this definition says that, if we are given two input/output pairs for which the inputs converge, and if the limits of the two inputs are close to each other, then the outputs become asymptotically close to each other. If $\mathcal{R}$ has an incremental limit gain $\kappa$, and if in addition $\mathcal{R}$ also admits a Cauchy gain, then both $\eta_1^\infty$ and $\eta_2^\infty$ exist whenever $\|\omega_1\|_{\mathrm{aa}} = \|\omega_2\|_{\mathrm{aa}} = 0$ (converging-input converging-output), and thus the "limsup" in Definition 1.4 is a limit, and the estimate becomes:

$$|\eta_1^\infty - \eta_2^\infty| \;\leq\; \kappa(|\omega_1^\infty - \omega_2^\infty|)\,. \tag{2}$$

With this concept, we have another obvious observation:

**Lemma 1.5** *(Small gain lemma for asymptotic amplitude, with uniqueness.)* Suppose that $\mathcal{R}$ and $\mathcal{S}$ are two behaviors with Cauchy gains $\gamma_1$ and $\gamma_2$ respectively, and incremental limit gains $\kappa_1$ and $\kappa_2$ respectively, and that the following condition holds:

$$\kappa_1(\kappa_2(r)) < r \quad \forall\, r > 0 \tag{3}$$

in addition to (1). Then, there exist two elements $\bar{u} \in \mathcal{U}$ and $\bar{y} \in \mathcal{Y}$ such that, for every input/output pair $(\omega, \eta) \in \mathcal{R} \bigcap \mathcal{S}^{-1}$ for which $\|\omega\|_{\text{aa}} < \infty$, $\omega^\infty = \bar{u}$ and $\eta^\infty = \bar{y}$.

*Proof.* If $\mathcal{R} \bigcap \mathcal{S}^{-1} = \emptyset$, there is nothing to prove. Otherwise, pick an arbitrary $(\omega_1, \eta_1) \in \mathcal{R} \bigcap \mathcal{S}^{-1}$ for which $\|\omega_1\|_{\text{aa}} < \infty$. From Lemma 1.2, there exist $\bar{u} := \omega_1^\infty$ and $\bar{y} := \eta_1^\infty$. Pick now any other $(\omega_2, \eta_2) \in \mathcal{R} \bigcap \mathcal{S}^{-1}$ for which $\|\omega_2\|_{\text{aa}} < \infty$; again by the Lemma, $\omega_2^\infty$ and $\eta_2^\infty$ exist. By the incremental limit gain property, in the form (2), both $|\bar{y} - \eta_2^\infty| \leq \kappa_1(|\bar{u} - \omega_2^\infty|)$ and $|\bar{u} - \omega_2^\infty| \leq \kappa_2(|\bar{y} - \eta_2^\infty|)$. From

$$|\bar{y} - \eta_2^\infty| \leq \kappa_1(\kappa_2(|\bar{y} - \eta_2^\infty|))$$

we conclude that $\eta_2^\infty = \bar{y}$, and so also $\omega_2^\infty = \bar{u}$. ∎

Once the appropriate definitions have been given, the two Lemmas are quite obvious. The harder step is, often, to verify when the Lemmas apply. In order to carry out such an application, one needs to find sufficient and easy to check conditions which guarantee the existence of Cauchy and incremental limit gains, for the systems whose feedback interconnection is being studied.

We will mainly study behaviors $\mathcal{R}$ which can be built up from cascades of simpler behaviors $\mathcal{R}_i$, each of which is either defined by a system of differential equations, by a pure delay, or by a memoryless nonlinearity. The composition $\mathcal{R}$ will represent the input/output pairs of a large set of delay-differential equations. The Cauchy and incremental limit gains of the behaviors $\mathcal{R}_i$ can be composed, so as to provide the gains of the complete system $\mathcal{R}$. We turn to that topic next.

## 2 Simple Behaviors and Cascades

The *delay-$\tau$ operator* $\mathcal{D}_\tau$ on $\mathcal{U}$, where $\tau \geq 0$, is the behavior, with $\mathcal{Y} = \mathcal{U}$, defined by: $(\omega, \eta) \in \mathcal{D}_\tau$ if and only if $\eta(t) = \omega(t - \tau)$ for all $t \geq \tau$. (The value of the output for $t < \tau$ is arbitrary; in an abstract dynamical systems sense, it forms part of the specification of initial conditions.) It is clear that $\mathcal{D}_\tau$ has Cauchy gain $I$ and incremental limit gain $I$, where $I$ is the identity function, $I(r) = r$.

Given any map $\psi : \mathcal{U} \to \mathcal{Y}$, the *memoryless behavior associated to $\psi$*, which we denote by $\mathcal{M}_\psi$, is the behavior consisting of all pairs of functions $(\omega, \eta)$ such that $\eta(t) = \psi(\omega(t))$ for all $t$. Suppose that $\psi$ is a Lipschitz map: for some $\lambda \geq 0$, $|\psi(u_1) - \psi(u_2)| \leq \lambda |u_1 - u_2|$ for all $u_1, u_2 \in \mathcal{U}$. Then $\mathcal{M}_\psi$ has Cauchy gain $\lambda I$ and incremental limit gain $\lambda I$, where $\lambda I(r) = \lambda r$.

Suppose that $\mathcal{R} \subseteq [\mathbb{R}_{\geq 0} \to \mathcal{U}] \times [\mathbb{R}_{\geq 0} \to \mathcal{Y}]$ and $\mathcal{S} \subseteq [\mathbb{R}_{\geq 0} \to \mathcal{Y}] \times [\mathbb{R}_{\geq 0} \to \mathcal{Z}]$ are two behaviors, with Cauchy gains $\gamma_1$ and $\gamma_2$ respectively, and consider the cascade combination shown pictorially in Figure 2 and defined formally as:

$$\mathcal{S} \circ \mathcal{R} := \{(\omega, \zeta) \mid (\exists\, \eta \in [\mathbb{R}_{\geq 0} \to \mathcal{Y}])\ \text{s.t.}\ (\omega, \eta) \in \mathcal{R}\ \&\ (\eta, \zeta) \in \mathcal{S}\}\,.$$

Then, clearly, $\mathcal{S} \circ \mathcal{R}$ has Cauchy gain $\gamma_2 \circ \gamma_1$. Suppose now that also $\mathcal{R}$ and $\mathcal{S}$ have incremental limit gains $\kappa_1$ and $\kappa_2$ respectively. Let $(\omega_i, \eta_i) \in \mathcal{R}$ and $(\eta_i, \zeta_i) \in \mathcal{S}$, $\|\omega_i\|_{\text{aa}} = 0$, for $i = 1, 2$. We have that $\eta_1^\infty$ and $\eta_2^\infty$ exist, and (2) holds with $\kappa = \kappa_1$. Similarly, since $\mathcal{S}$ has a Cauchy gain, $\zeta_1^\infty$ and $\zeta_2^\infty$ exist, and $|\zeta_1^\infty - \zeta_2^\infty| \leq \kappa_2(|\eta_1^\infty - \eta_2^\infty|)$. Therefore

$$|\zeta_1^\infty - \zeta_2^\infty| \leq \kappa_2(\kappa_1(|\omega_1^\infty - \omega_2^\infty|))$$

and hence $\mathcal{S} \circ \mathcal{R}$ has incremental limit gain $\kappa_2 \circ \kappa_1$.

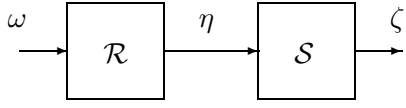

Figure 2: Cascade $\mathcal{S} \circ \mathcal{R}$

## 2.1 Tighter Estimates: Relative Gains

Tighter estimates of gains for the cascade $\mathcal{S} \circ \mathcal{R}$ can make use of the following observation. Suppose that the possible output signals of $\mathcal{R}$ all tend, as $t \to \infty$, to values in a restricted subset $Y$ of $\mathcal{Y}$. Then the relevant gains $\gamma_2, \kappa_2$ should be the gains of $\mathcal{S}$ when restricted to those signals of the form $(\eta, \zeta) \in \mathcal{S}$ such that $\eta \in [\mathbb{R}_{\geq 0} \to Y]$. These gains may well be smaller that the original ones, so that smaller overall gains result for the cascade. Let us make this precise.

For any subset $U_0 \subseteq \mathcal{U}$, we write "$\omega \to U_0$" if $\omega(t)$ converges to $U_0$ as $t \to \infty$, that is, for every $\varepsilon > 0$ there is some $T \geq 0$ such that

$$\omega(t) \in B_\varepsilon(U_0) = \{u \in \mathcal{U} \mid (\exists\, u' \in U_0) \, |u - u'| \leq \varepsilon\}$$

for every $t \geq T$.

Let $U_0 \subseteq \mathcal{U}$ and let $\mathcal{R} \subseteq [\mathbb{R}_{\geq 0} \to \mathcal{U}] \times [\mathbb{R}_{\geq 0} \to \mathcal{Y}]$. We will say that $\mathcal{R}$ has a *Cauchy gain* $\gamma$ *on* $U_0$ if $\|\eta\|_{\text{aa}} \leq \gamma(\|\omega\|_{\text{aa}})$ holds for each input/output pair $(\omega, \eta) \in \mathcal{R}$ for which $\omega \to U_0$. Similarly, we say that $\mathcal{R}$ has *incremental limit gain* $\kappa$ *on* $U_0$ if $\limsup_{t\to\infty} |\eta_1(t) - \eta_2(t)| \leq \kappa(|\omega_1^\infty - \omega_2^\infty|)$ holds whenever $(\omega_i, \eta_i) \in \mathcal{R}$ are any two pairs such that $\omega_1^\infty$ and $\omega_2^\infty$ both exist and belong to $U_0$. In the special case $U_0 = \mathcal{U}$, one recovers the definitions of Cauchy and incremental limit gains.

Suppose now that there are two sets $U_0 \subseteq \mathcal{U}$ and $Y_0 \subseteq \mathcal{Y}$ such that:

- $\mathcal{R}$ has Cauchy gain $\gamma_1$ on $U_0$.
- $\mathcal{S}$ has Cauchy gain $\gamma_2$ on $Y_0$.
- Whenever $(\omega, \eta) \in \mathcal{R}$ is so that $\omega \to U_0$, necessarily $\eta \to Y_0$.

Then, clearly, $\mathcal{S} \circ \mathcal{R}$ has Cauchy gain $\gamma_2 \circ \gamma_1$ on $U_0$. An analogous conclusion holds for incremental limit gain on $U_0$.

## 2.2 A Sufficient Condition

Recall that, for any metric space $\mathcal{U}$ and function $\omega : \mathbb{R}_{\geq 0} \to \mathcal{U}$, the omega-limit set $\Omega = \Omega^+[\omega]$ is the set consisting of those points $u \in \mathcal{U}$ for which there exists a convergent sequence $\omega(t_i) \to u$, for some sequence $\{t_i\} \subseteq \mathbb{R}_{\geq 0}$ such that $t_i \to \infty$ as $i \to \infty$. The following properties are elementary: (1) the set $\Omega$ is closed; (2) if $\omega \to U$ and $U$ is closed, then $\Omega \subseteq U$; and (3) if $\omega$ is precompact, that is to say, if there is some compact subset $U \subseteq \mathcal{U}$ such that $\omega(t) \in U$ for all $t \geq 0$, then $\Omega$ is compact, and $\omega \to \Omega$ (proof of this last statement: if there is some $\varepsilon > 0$ and some sequence $t_i \to \infty$ such that $\omega(t_i) \in U \setminus B_\varepsilon(\Omega)$ for all $i$, then one can pick a subsequence of $\{t_i\}$ such that $\omega(t_{i_j}) \to u$ for some $u$, and thus $u \in U \setminus B_{\varepsilon/2}(\Omega)$, a contradiction since $u \in \Omega$ by definition of $\Omega$).

In general, we denote by $|U|$ the *diameter* $\sup\{|u - v| \mid u, v \in U\}$ of a closed subset $U$ of a metric space $\mathcal{U}$. For each $\omega : \mathbb{R}_{\geq 0} \to \mathcal{U}$, it holds that $|\Omega^+[\omega]| \leq \|\omega\|_{\text{aa}}$, and, if $\omega$ is precompact,

$$\|\omega\|_{\text{aa}} \;=\; |\Omega^+[\omega]| \,. \tag{4}$$

Indeed, pick any $\varepsilon > 0$ and two elements $u, v \in \Omega$ such that $|u - v| \geq |\Omega| - \varepsilon$; then there are two sequences $\omega(t_i) \to u$ and $\omega(s_i) \to v$, so $\|\omega\|_{\text{aa}} = \limsup_{s,t \to \infty} |\omega(t) - \omega(s)| \geq |u - v| \geq |\Omega| - \varepsilon$. As this is true for every $\varepsilon > 0$, we have $\|\omega\|_{\text{aa}} \geq |\Omega|$. Conversely, if $|\omega(t_i) - \omega(s_i)| \geq \|\omega\|_{\text{aa}} - \varepsilon$ for some two sequences $t_i \to \infty$ and $s_i \to \infty$, we may extract first a subsequence of $\{t_i\}$ such that $\omega(t_{i_j})$ is convergent (precompactness is used here), and then a subsequence of $\{s_{i_j}\}$, so that, without loss of generality we may suppose that $\omega(t_i) \to u$ and $\omega(s_i) \to v$ for some $u, v \in \Omega$, and thus $|\Omega| \geq |u - v| \geq \|\omega\|_{\text{aa}} - \varepsilon$, so letting $\varepsilon \to 0$ gives the other inequality.

We say that a mapping $\Gamma$ assigning subsets of one set to subsets of another is monotonic if $U_1 \subseteq U_2 \Rightarrow \Gamma(U_1) \subseteq \Gamma(U_2)$.

**Lemma 2.1** Suppose given a behavior $\mathcal{R}$, a compact subset $U_0 \subseteq \mathcal{U}$, a function $\gamma \in \mathcal{K}_\infty$, and a mapping $\Gamma$ from compact subsets of $U_0$ to subsets of $\mathcal{Y}$, such that the following properties hold:

(a) For each $(\omega, \eta) \in \mathcal{R}$ for which $\Omega^+[\omega] \subseteq U_0$, the output $\eta$ is precompact.

(b) For each compact subset $U \subseteq U_0$, and each $(\omega, \eta) \in \mathcal{R}$ for which $\Omega^+[\omega] \subseteq U$, it holds that $\Omega^+[\eta] \subseteq \Gamma(U)$.

(c) For each compact subset $U \subseteq U_0$, it holds that $|\Gamma(U)| \leq \gamma(|U|)$.

Then $\mathcal{R}$ has Cauchy gain $\gamma$ on $U_0$ and incremental limit gain $\gamma$ on $U_0$. Moreover, for each compact subset $U \subseteq U_0$, and each $(\omega, \eta) \in \mathcal{R}$ for which $\omega \to U$, $\eta \to \Gamma(\Omega^+[\omega])$. If $\Gamma$ is monotonic, then also $\eta \to \Gamma(U)$.

*Proof.* Pick any $(\omega, \eta) \in \mathcal{R}$ and any compact $U \subseteq U_0$, and suppose that $\omega \to U$. By (1) and (2) in the previous discussion, the set $\Omega^+[\omega]$ is a compact (since closed) subset of $U$. By (a), $\eta$ is precompact. Therefore $\|\eta\|_{\text{aa}} = |\Omega^+[\eta]|$, and also $\eta \to \Omega^+[\eta]$. By (b), applied to $\Omega^+[\omega]$ itself, we know that $\Omega^+[\eta] \subseteq \Gamma(\Omega^+[\omega])$, which gives the conclusion $\eta \to \Gamma(\Omega^+[\omega])$. If $\Gamma$ is monotonic, then $\Omega^+[\omega] \subseteq U$ implies that $\Gamma(\Omega^+[\omega]) \subseteq \Gamma(U)$, so $\eta \to \Gamma(U)$. In addition, $|\Omega^+[\eta]| \leq |\Gamma(\Omega^+[\omega])|$ together with (c) give the following inequality:

$$\|\eta\|_{\text{aa}} \;=\; |\Omega^+[\eta]| \;\leq\; |\Gamma(\Omega^+[\omega])| \;\leq\; \gamma(|\Omega^+[\omega]|) \;\leq\; \gamma(\|\omega\|_{\text{aa}}) \,.$$

When applied in the special case $U = U_0$, this establishes the Cauchy gain conclusion.

Suppose now that $(\omega_i, \eta_i) \in \mathcal{R}$ are any two pairs such that $\omega_1^\infty$ and $\omega_2^\infty$ both exist and belong to $U_0$. In particular, $\omega_1 \to U_0$ and $\omega_2 \to U_0$. So both $\eta_1^\infty$ and $\eta_2^\infty$ exist, by the Cauchy gain conclusion. Note that $\Omega^+[\omega_i] = \{\omega_i^\infty\}$ and $\Omega^+[\eta_i] = \{\eta_i^\infty\}$ for $i = 1, 2$. We introduce the two-element set $U = \{\omega_1^\infty, \omega_2^\infty\} \subseteq U_0$; note that $|U| = |\omega_1^\infty - \omega_2^\infty|$. From $\Omega^+[\omega_i] \subseteq U$ and (b), we have that $\Omega^+[\eta_i] \subseteq \Gamma(U)$, that is, $\eta_i^\infty \in \Gamma(U)$, for $i = 1, 2$. Therefore

$$|\eta_1^\infty - \eta_2^\infty| \leq |\Gamma(U)| \leq \gamma(|U|) = \gamma(|\omega_1^\infty - \omega_2^\infty|),$$

which proves the incremental limit property. ∎

## 3 Systems of Differential Equations

A particular class of behaviors, in fact the main objects of interest in this note, are obtained as follows. We consider systems of differential equations with inputs and outputs:

$$\dot{x} = f(x, u), \quad y = h(x) \tag{5}$$

for which states $x(t)$ evolve in a subset $\mathcal{X}$ of a Euclidean space $\mathbb{R}^n$, inputs take values $u(t)$ in a complete metric space $\mathcal{U}$ and outputs take values $y(t)$ in a complete metric space $\mathcal{Y}$. (Typically in applications, $\mathcal{U}$ and $\mathcal{Y}$ are any two closed subsets of Euclidean spaces.) Technically, we assume that $f : \mathcal{X}_0 \times \mathcal{U} \to \mathbb{R}^n$ is defined on an open subset $\mathcal{X}_0 \subseteq \mathbb{R}^n$ which contains $\mathcal{X}$, is continuous, and is locally Lipschitz in $x$ uniformly on compact subsets of $\mathcal{X}_0 \times \mathcal{U}$; the map $h : \mathcal{X} \to \mathcal{Y}$ is assumed to be continuous. Furthermore, $\mathcal{X}$ is an invariant and forward complete subset, in the sense that, for each Lebesgue-measurable precompact input $\omega : \mathbb{R}_{\geq 0} \to \mathcal{U}$, and each initial state $x_0 \in \mathcal{X}$, the unique solution $\xi(t) = \varphi(t, x_0, \omega)$ of the initial value problem $\dot{\xi}(t) = f(\xi(t), \omega(t))$, $\xi(0) = x_0$, is defined and satisfies $\xi(t) \in \mathcal{X}$ for all $t \geq 0$. (The function $\xi$ is Lipschitz, and hence differentiable almost everywhere; if $\omega$ is continuous, then $\xi$ is continuously differentiable.) To any given system (5) one associates a behavior $\mathcal{R}$, with input-value space $\mathcal{U}$ and output-value space $\mathcal{Y}$, defined by: $(\omega, \eta) \in \mathcal{R}$ if and only if $\omega$ is precompact and Lebesgue-measurable, and there exists some $x_0 \in \mathcal{X}$ such that $\eta(t) = h(\varphi(t, x_0, \omega))$ for all $t \in \mathbb{R}_{\geq 0}$. We call $\mathcal{R}$ *the behavior of* (5).

**Remark 3.1** A minor technicality concerns the fact that Lebesgue-measurable functions are, strictly speaking, not functions but equivalence classes of functions, so one should interpret the "limsup" in the definition of asymptotic amplitude in an "almost everywhere" manner; similarly, we interpret "precompact" as meaning that there is some $\omega$ in the given equivalence class whose values all remain in a compact. From now on, we leave this technicality implicit; in applications to stability of feedback loops involving systems of differential equations, all the functions considered are continuous –even differentiable– so the issue does not even arise. □

We will obtain sufficient conditions for the existence of the two types of gains, expressed in terms of Lyapunov-type functions.

Given a subset $U \subseteq \mathcal{U}$, we will say that a function

$$V : \mathcal{X} \to \mathbb{R}_{\geq 0}$$

is a *U-decrease* function provided that the following properties hold:

- $V$ is proper, that is, $\{x \in \mathcal{X} \mid a \leq V(x) \leq b\}$ is a compact subset of $\mathcal{X}$, for each $a \leq b$ in $\mathbb{R}_{\geq 0}$;

- $V$ is continuous;

- for each $x \in \mathcal{X}$ which does not belong to $Z_V := \{x \mid V(x) = 0\}$, it holds that $V$ is continuously differentiable in a neighborhood of $x$ and

$$\nabla V(x) \cdot f(x, u) < 0 \quad (6)$$

for all $u \in U$.

(We understand continuous differentiability in the following sense: there is a neighborhood of $x$ in $\mathcal{X}_0$ such that $V$ extends to this neighborhood as a $\mathcal{C}^1$ function.)

**Lemma 3.2** Suppose that $V$ is a $U$-decrease function, for some compact subset $U \subseteq \mathcal{U}$. Pick any Lebesgue-measurable precompact $\omega : \mathbb{R}_{\geq 0} \to \mathcal{U}$ and any solution $\xi$ of the system $\dot{\xi} = f(\xi, \omega)$. Suppose that either:

1. there is some $T_0 \geq 0$ such that $\omega(t) \in U$ for all $t \geq T_0$, or

2. $\omega \to U$ and $\xi$ is precompact.

Then $\xi \to Z_V$.

*Proof.* Given any $\omega$ and $\xi$, we will first let $a > 0$ be arbitrary and prove that the set $V_a := \{x \mid V(x) \leq a\}$ has the property that, for some $T^* \geq 0$,

$$\xi(t) \in V_a \quad \forall\, t \geq T^*. \quad (7)$$

If the assumption is that $\omega \to U$ as $t \to \infty$ and that $\xi$ is precompact, that is to say, there is some compact subset $C_0$ of $\mathcal{X}$ such that $\xi(t) \in C_0$ for all $t \geq 0$, then we introduce $b := \max\{V(x), x \in C_0\}$ and the set $C := V^{-1}([a, b])$. Note that $\xi(t) \in C$ whenever $V(\xi(t)) \geq a$, by the choice of $b$. Since, by properness of $V$, $C$ is a compact subset of $\mathcal{X} \setminus Z_V$, by Property (6) there is some $\alpha > 0$ and some neighborhood $\widetilde{U}$ of $U$ in $\mathcal{U}$ so that $\nabla V(x) \cdot f(x, u) \leq -\alpha$ for all $x \in C$ and all $u \in \widetilde{U}$. Since $\omega \to U$, there must be some $T_0$ such that $\omega(t) \in \widetilde{U}$ for all $t \geq T_0$. Thus,

$$\nabla V(x) \cdot f(x, \omega(t)) \leq -\alpha < 0 \quad \forall\, x \in C \quad \forall\, t \geq T_0. \quad (8)$$

If, instead, the assumption is that there is some $T_0 \geq 0$ such that $\omega(t) \in U$ for all $t \geq T_0$, we pick $T_0$ to be any such number, and let $C := \{x \mid V(x) \geq a\}$. So once more we have that $\nabla V(x) \cdot f(x, \omega(t)) < 0$ for all $x \in C$ and all $t \geq T_0$.

*Claim:* Let $T_1 := \inf\{t \geq T_0 \mid V(\xi(t)) \leq a\}$ (if $V(\xi(t)) > a$ for all $t \geq T_0$, we define $T_1 = +\infty$). Then $V(\xi(t)) \leq V(\xi(T_0))$ for all $t \in [T_0, T_1]$ and $V(\xi(t)) \leq a$ for all $t \geq T_1$.

*Proof of the claim:* Suppose that $\xi(t) \in C$ for all $t$ in some interval $(\tau_1, \tau_2)$ with $\tau_1 \geq T_0$. Then $dV(\xi(t))/dt = \nabla V(\xi(t)) \cdot f(\xi(t), \omega(t)) < 0$ for almost every $t \in (\tau_1, \tau_2)$. Therefore, $V(\xi(t))$ is decreasing on this interval, and we have that $V(\xi(t)) \leq V(\xi(\tau_1))$ for all $t \in (\tau_1, \tau_2)$. In particular, for each $t \in (T_0, T_1)$, by minimality of $T_1$ we know that $V(\xi(t)) > a$ and so $\xi(t) \in C$. This proves the first part of the claim: $V(\xi(t)) \leq V(\xi(T_0))$ for all $t \in [T_0, T_1]$. If $T_1 = \infty$, there is nothing more to prove. So assume that $T_1 < \infty$ and there exists some $S > T_1$ such that $V(\xi(S)) > a$. Then there is some $T \in [T_1, S)$ such that $V(\xi(T)) = a$. We pick $T' \in [T_1, S)$ to be maximal with this property. It follows that $V(\xi(t)) > a$ for all $t \in (T', S]$. Applying the above argument with $\tau_1 = T'$ and $\tau_2 = S$, we have that $V(\xi(S)) \leq a = V(\xi(T'))$, a contradiction. So the claim holds.

We conclude that $V(\xi(t)) \leq \max\{a, V(\xi(0))\}$ for all $t \geq T_0$. Therefore the trajectory $\xi$ is precompact, and the first case in the Lemma is included in the second case, so we can assume that (8) holds. We claim that this means that $T_1 < \infty$, so that (7) holds with $T^* = T_1$. Indeed, if this were not true, then $\xi(t) \in C$ for all $t \geq T_0$, so $dV(\xi(t))/dt \leq -\alpha$ for almost all $t$, which gives $V(\xi(t)) \leq V(\xi(0)) - \alpha t$ for all $t \geq T_0$, which is impossible since $V$ is nonnegative.

To conclude that $\xi \to Z_V$, since $\xi$ is precompact we need only show that its omega-limit set $\Omega^+[\xi]$ is contained in $Z_V$. To see this, we pick any $z \in \Omega^+[\xi]$ and a sequence $\xi(t_i) \to z$. So $V(\xi(t_i)) \to V(z)$. If $z \notin Z_V$, let $a := V(z)/2 \neq 0$. Then Property (7) gives that $V(\xi(t_i)) \leq a$ for all $i$ large enough, which says that $\limsup_i V(\xi(t_i)) \leq a$, contradicting $V(\xi(t_i)) \to 2a$. Thus $z \in Z_V$. ∎

**Theorem 1** *Suppose that $\mathcal{R}$ is the behavior of a system (5), and there is some $\gamma \in \mathcal{K}_\infty$ such that the following property holds: for each compact subset $U \subseteq \mathcal{U}$, there exists a $U$-decrease function $V_U$ such that*
$$|h(Z_{V_U})| \leq \gamma(|U|).$$
*Then, for each compact subset $U \subseteq \mathcal{U}$,*

1. *if $(\omega, \eta) \in \mathcal{R}$ is such that $\omega \to U$, then $\eta \to h(Z_{V_U})$;*

2. *$\mathcal{R}$ has Cauchy gain $\gamma$ on $U$ and incremental limit gain $\gamma$ on $U$.*

*In particular, $\mathcal{R}$ has Cauchy gain $\gamma$ and incremental limit gain $\gamma$.*

*Proof.* We will first show that, for every triple $(\omega, \xi, \eta)$ with $\omega$ precompact and Lebesgue-measurable, such that $\dot\xi(t) = f(\xi(t), \omega(t))$ and $\eta(t) = h(\xi(t))$, the following properties hold:

(i) $\xi$ and $\eta$ are precompact;

(ii) for each compact $U \subseteq \mathcal{U}$ and each $U$-decrease function $V$, if $\omega \to U$ then $\xi \to Z_V$ and $\eta \to h(Z_V)$.

Since $\omega$ is precompact, there exists some compact set, let us call it $U'$, such that $\omega(t) \in U'$ for all $t$. As there exists some $U'$-decrease function $V'$, the first case in Lemma 3.2, applied to this data and with $T_0 = 0$, gives that $\xi(t) \to Z_{V'}$, so, being continuous as a function of $t$, $\xi$ is precompact. Next we apply once more Lemma 3.2, using now the second case with any given compact $U$ and a $U$-decrease function $V$, to conclude that $\xi \to Z_V$. Since the set $Z_V$ is compact and the mapping $h$ is continuous, it follows that also $\eta(t) = h(\xi(t)) \to h(Z_V)$ as $t \to \infty$, and $\eta$ is precompact as well.

Now we pick any compact $U_0 \subseteq \mathcal{U}$, and for each compact subset $U \subseteq U_0$ we pick a $U$-decrease function $V$ such that $|h(Z_V)| \leq \gamma(|U|)$, and let $\Gamma(V) := h(Z_V)$. We will apply Lemma 2.1. Property (c) in that Lemma holds by definition of $\Gamma$. Also, (a) holds, by (i). To prove that (b) is true, suppose that $\Omega^+[\omega] \subseteq U$. Since $\omega$ is precompact, $\omega \to \Omega^+[\omega]$, so also $\omega \to U$. By (ii), we know that $\eta \to h(Z_V)$, and so $\Omega^+[\eta] \subseteq h(Z_V) = \Gamma(V)$, as desired. The Lemma then says that $\mathcal{R}$ has Cauchy gain $\gamma$ on $U_0$ and also has incremental limit gain $\gamma$ on $U_0$.

Finally, given an arbitrary pair $(\omega, \eta) \in \mathcal{R}$, we pick some compact $U_0$ such that $\omega(t) \in U_0$; then the Cauchy property on $U_0$ gives that $\|\eta\|_{\mathrm{aa}} \leq \gamma(\|\omega\|_{\mathrm{aa}})$, and similarly for the incremental gain. ∎

## 4 An Example

As an illustration of gain computations and small-gain stability arguments, we consider systems which consist of cascades of several subsystems, each of which can be individually described by some ordinary differential equation $\dot{x}_i = f(x_i, u_i)$ with input $u_i$. The input $u_1 = u$ to the first of the systems in the cascade is an external one, while the intermediate inputs $u_i$, $i > 1$, between two stages depend on the state of the preceding stage. For instance, $x_i(t)$ might represent the amount present, at any given time $t$, of the activated form $E_i^*$ of an enzyme $E_i$ whose production rate is, in turn, dependent upon the amount present of the activated form $E_{i-1}^*$ of the enzyme $E_{i-1}$. We allow transport delays in between stages. This leads to systems given by sets of equations as follows:

$$\begin{aligned}
\dot{x}_1(t) &= f_1(x_1(t), u(t)) \\
\dot{x}_2(t) &= f_2(x_2(t), x_1(t - \tau_1)) \\
&\vdots \\
\dot{x}_n(t) &= f_n(x_n(t), x_{n-1}(t - \tau_{n-1}))
\end{aligned}$$

where $\tau_1, \ldots, \tau_{n-1} \geq 0$ are the delays among the stages of the process (the particular case in which there are no delays is included in this formalism by setting all $\tau_i = 0$). See Figure 3, where $x_i^\tau(t) := x_i(t - \tau)$ and we use $\mathcal{R}_i$ to denote the behavior associated to the system $\dot{x} = f_i(x, u)$ with output $y = x$. One often asks about such systems whether adding a feedback loop from the last stage to the first, as shown in Figure 4, might introduce instabilities, such as oscillations or even chaotic behavior. Specifically, one may have, for instance, that the action of $u$ on the first subsystem is inhibited by the final product $x_n$. Assuming that all the variables $x_i$, as well as the external input $u$, are scalar and take only nonnegative values (as is the case when they represent concentrations of an enzyme), a typical model for inhibition is obtained when the first equation becomes

$$\dot{x}_1(t) \;=\; f\left(x_1(t), \frac{u(t)}{1+kx_n(t-\tau_n)}\right)$$

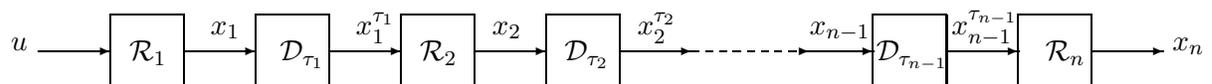

Figure 3: Cascade of $\mathcal{R}_i$'s and Delays

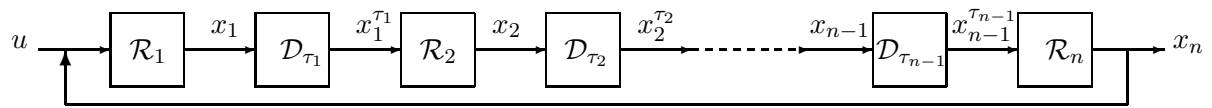

Figure 4: Cascade of $\mathcal{R}_i$'s and Delays, With Feedback to First Stage

(other expressions for inhibition are also possible, of course), where $k \geq 0$ serves to parametrize the feedback strength. Note that we are also allowing for an additional delay in the feedback. Suppose that we are interested in analyzing the case in which $u(t)$ equals a constant value $\mu$, so that the effective input being fed to the first subsystem is $\omega(t) = \psi(x_n^{\tau_n}(t)) = \psi(x_n(t-\tau_n))$, where $\psi$ is the function:

$$\psi(x) \;=\; \frac{\mu}{1+kx}\,.$$

The closed-loop system may be viewed as the feedback interconnection of the memoryless behavior $\mathcal{M}_\psi$ with the behavior of the forward composite system with output $x_n^{\tau_n}$, that is, the composition

$$\mathcal{R} \;=\; \mathcal{D}_{\tau_n} \circ \mathcal{R}_n \circ \mathcal{D}_{\tau_{n-1}} \circ \ldots \circ \mathcal{D}_{\tau_1} \circ \mathcal{R}_1 \,.$$

Suppose that each $\mathcal{R}_i$ has Cauchy gain $\gamma_i$ and also incremental limit gain $\gamma_i$. Then $\mathcal{R}$ has gain $\gamma = \gamma_n \circ \ldots \circ \gamma_1$ of both types. Since $\psi$ has Lipschitz constant $k\mu$, $\mathcal{M}_\psi$ has both gains $k\mu I$. Therefore, provided that the small gain condition

$$\gamma(k\mu r) < r \quad \forall\, r > 0$$

holds, one concludes from Lemma 1.5 that there is some value $\bar{u}$ such that, for every solution of the closed-loop system, $\omega \to \bar{u}$. This means, in turn, because each $\mathcal{R}_i$ has a Cauchy gain, that every state variable $x_i$ converges to a unique equilibrium (independently of initial conditions).

The problem, therefore, is to estimate gains $\gamma_i$ for the systems $\mathcal{R}_i$. We briefly discuss one situation, itself of great interest, in which estimates can be obtained.

Suppose given intervals $\mathcal{X}_i = [a_i, b_i]$ and $\mathcal{U}_i = [u_i, v_i]$, $i = 1, \ldots, n$, such that any solution of $\dot{x} = f_i(x, u)$ with initial conditions in $\mathcal{X}_i$ and input with values in $\mathcal{U}_i$ remains in $\mathcal{X}_i$. For example, let $\mathcal{X}_i = [0, x_i^{\max}]$, where $x_i^{\max}$ is the maximum possible amount of a substance, such as the activated form of an enzyme, that may be synthesized. Also, suppose given, for each $i$, a strictly increasing and onto function

$$g_i \,:\, [a_i, b_i) \to \mathbb{R}_{\geq 0}$$

with the following properties:

1. each $g_i^{-1} : [0, \infty) \to [a_i, b_i)$ is Lipschitz, with Lipschitz constant $\lambda_i$;

2. $x < g_i^{-1}(u) \Rightarrow f_i(x, u) > 0$, and $x > g_i^{-1}(u) \Rightarrow f_i(x, u) < 0$, for every $u \in \mathcal{U}_i$ and $x \in \mathcal{X}_i$.

For any given interval $U = [c, d] \subseteq \mathcal{U}_i$, we may introduce the function $V$ which measures the distance from any $x$ to the set $g^{-1}(U) = [g_i^{-1}(c), g_i^{-1}(d)]$. Thus $V$ is a $U$-decrease function, and $Z_V = [g_i^{-1}(c), g_i^{-1}(d)]$ has diameter $|Z_V| \leq \lambda\, |d - c| = \lambda\, |U|$. We conclude that $\mathcal{R}$ admits the Cauchy gain $\gamma(r) = \lambda_1 \ldots \lambda_n r$, and also incremental limit gain $\gamma$. Hence the desired stability result holds as long as

$$k \;<\; \frac{1}{\mu\,\lambda_1 \ldots \lambda_n}$$

is satisfied by the gain $k$.

In particular, assume that each $f_i$ has the following form, for different $\alpha$'s and $\beta$'s: $f(x, u) = -\alpha(x) + u\beta(x)$, where $\alpha$ and $\beta$ are nonnegative on $\mathcal{X} = [a, b]$, $\alpha$ is strictly increasing and $\beta$ is strictly decreasing, and $\alpha(a) = \beta(b) = 0$. These conditions guarantee that $[a, b]$ is invariant when $u(t) \in \mathbb{R}_{\geq 0}$, and the function $g(x) := \frac{\alpha(x)}{\beta(x)}$ is as required.

Note that the effective input signal $\omega$ is not arbitrary, but is restricted to lie in the interval $U_0 = [\frac{\mu}{(1+kb_n)}, \frac{\mu}{(1+ka_n)}]$. Thus $x_1$ will approach the set $U_1 = g^{-1}(U_1)$, and so forth inductively. This allows the use of relative gain results, and far tighter estimates on gains, and will be developed in detail for an example, in the extended version of this report.